\documentclass[a4paper,11pt]{article}

\usepackage[cmex10]{amsmath}
\usepackage{amsfonts}
\usepackage{amssymb}
\usepackage{amsthm}
\usepackage[mathscr]{euscript}

\usepackage[all,2cell,ps]{xy}   

\usepackage{url}	  

\usepackage[dvips]{graphicx}    

\theoremstyle{definition} \newtheorem{theorem}{Theorem}[section]
\theoremstyle{definition} \newtheorem{definition}[theorem]{Definition}
\theoremstyle{definition} \newtheorem{lemma}[theorem]{Lemma}
\theoremstyle{definition} \newtheorem{proposition}[theorem]{Proposition}
\theoremstyle{definition} 
\theoremstyle{definition} 
\theoremstyle{definition} 
\theoremstyle{definition} 
\theoremstyle{definition} \newtheorem{remark}[theorem]{Remark}
\theoremstyle{definition} 
\newcommand{\eat}[1]{}

\def\pf{\begin{proof}}

\def\sqr{\end{proof}}

\begin{document}

\title{\sc The Semiring Properties of\\ Boolean Propositional Algebras}
\date{}

\author{Mahesh Rudrachar\footnote{Unisys Global Services India,
    Bangalore 560 025, India, {\small \tt mahesh.rudrachar@in.unisys.com} } 
    \ \ \ 
    Shrisha Rao\footnote{International Institute of Information Technology -
    Bangalore, Bangalore 560 100, India, {\small \tt
      srao@iiitb.ac.in}} \ \ \ 
    Amit Raj\footnote{Unisys Global Services India,
    Bangalore 560 025, India, {\small \tt amit.raj@in.unisys.com} }}
\maketitle

\begin{abstract}

  This paper illustrates the relationship between boolean
  propositional algebra and semirings, presenting some results of
  partial ordering on boolean propositional algebras, and the
  necessary conditions to represent a boolean propositional subalgebra
  as equivalent to a corresponding boolean propositional algebra.  It
  is also shown that the images of a homomorphic function on a boolean
  propositional algebra have the relationship of boolean propositional
  algebra and its subalgebra.  The necessary and sufficient conditions
  for that homomorphic function to be onto-order preserving, and also
  an extension of boolean propositional algebra, are explored.
  
\end{abstract}


\section{Introduction}

The English mathematician George Boole (1815--1864) sought to give
symbolic form to Aristotle's system of logic. Boole wrote a treatise
on the subject in 1854, titled \emph{An Investigation of the Laws of
  Thought, on Which Are Founded the Mathematical Theories of Logic and
  Probabilities}, which codified several rules of relationship between
mathematical quantities limited to one of two possible values: true or
false, 1 or 0. His mathematical system became known as \emph{Boolean
  algebra}.

Bourne~\cite{Bourne} has discussed the homomorphism theorems for
semirings.  Allen~\cite{Allen} has discussed the extension of a
theorem of Hilbert to semirings.  Zeleznikow~\cite{Zelez} has
discussed the natural partial order on semirings.  This paper
illustrates the general idea of interpreting properties of boolean
propositional algebras, including partial ordering, homomorphism,
isomorphism and differences, when they are taken as semirings.

Let \(p, q, r\) be propositions, and $\mathcal{P}$ be the set of all
propositions in the universe of discourse. A \emph{Boolean
  Propositional Algebra (BPA)} $\mathcal{B}$ ~\cite{Eldon} is a
six-tuple consisting of a set $\mathcal{P}$, equipped with two binary
operations $\wedge$ (called `meet' or `and') and $\vee$ (called `join'
or `or'), a unary operation $\neg$ (called `complement' or `not') and
two elements $0$ and $1$ and it is denoted by
$(\mathcal{P},\wedge,\vee,\neg,0,1)$.  If we use infix operators like
$<$ or $\preccurlyeq$ to compare propositions, we may write
$<_{\mathcal{B}}$ and $\preccurlyeq_{\mathcal{B}}$ to clarify which
BPA is taken as the scope.

If there is some subset $\widehat{\mathcal{B}}$ of a boolean
propositional algebra $\mathcal{B}$ that is a BPA in its own right,
then it may be called a \emph{Boolean Propositional Subalgebra (BPSA)}
~\cite{Frank} If we use infix operators like $<$ or $\preccurlyeq$, we
may write $<_{\widehat{\mathcal{B}}}$ and
$\preccurlyeq_{\widehat{\mathcal{B}}}$ to indicate the corresponding
scope.


\section{Boolean Propositional Algebra of Monoids and
  Semirings} \label{monoid}

It is possible to define an arithmetic on propositional logic in the
obvious way: take $+$ to mean $\wedge$, and $\times$ to mean $\vee$.

\begin{proposition} \label{opsproperties}

Given propositions $p$, $q$, and $r$, we have the following.

\begin{itemize}

\item[(i)] $+$ and $\times$ are commutative and associative: \(p
  \times q = q \times p\); \(p + r = r + p\); \(p + (q + r) = (p + q)
  + r\); and \(p \times (q \times r) = (p \times q) \times r\).

\item[(ii)] $\times$ distributes over $+$: \(p \times (q + r) = (p
  \times q) + (p \times r)\).

\end{itemize}
\end{proposition}

Given this propositional arithmetic, we can posit the existence of two
identity operators, one each for $+$ and $\times$.

\begin{definition}

  The multiplicative and additive identities are defined as follows.

\begin{itemize}

\item[(i) ] The additive identity \(\top\) is the proposition such
  that for any proposition $p$, \(p + \top = \top + p = p\).

\item[(ii) ] The multiplicative identity \(\bot\) is the proposition
  such that for any proposition $p$, \(p \times \bot = \bot \times p =
  p\).

\end{itemize}

\end{definition}

By the commutativity of the $+$ and $\times$ operators, we observe
that the identity elements are two-sided.

Informally, we may describe these elements as follows:

\begin{itemize}

\item[(i)] The additive identity \(\top\) is a proposition ``that is
  always true.''  The direct sum of such a proposition and $p$ is
  obviously $p$ itself.

\item[(ii)] The multiplicative identity \(\bot\) is a proposition
  ``that is always false.''  The direct product of such a proposition
  and $p$ is likewise $p$ itself.

\end{itemize}

Then $\mathcal{P}$, combined with the $+$ operator, is a \emph{monoid}
(a set with an associative operator and a two-sided identity
element)~\cite{Hungerford}.  Similarly, $\mathcal{P}$ is also a monoid
when considering the $\times$ operator.  For notational convenience,
we denote these monoids as \((\mathcal{P},+)\) and
\((\mathcal{P},\times)\).


It is further clear that the set $(\mathcal{P},+,\times)$ is a
\emph{semiring} when taken with the operations $+$ and $\times$
because the following conditions~\cite{Golan1999} for being a semiring
are satisfied:

\begin{itemize}

\item[(i)] \((\mathcal{P}, +)\) is a commutative monoid with identity
  element \(\top\);

\item[(ii)] \((\mathcal{P}, \times)\) is a monoid with identity
  element \(\bot\);

\item[(iii)] \(\times\) distributes over \(+\) from either side;

\item[(iv)] \(\top \times p = \top = p \times \top\)
  for all \(p \in \mathcal{P}\).

\end{itemize}

This proposition semiring will be denoted by
\((\mathcal{P},+,\times)\), and its properties are as indicated in the
following.

\begin{remark} \label{remark1}

  The semiring $(\mathcal{P},+,\times)$ is \emph{zerosumfree}, because
  \(p + q = \top\) implies, for all \(p, q \in \mathcal{P}\), that \(p
  = q = \top\).

\end{remark}

This property shows~\cite{Golan1999} that the monoid
\((\mathcal{P},+)\) is completely removed from being a group, because
no non-trivial element in it has an inverse.

A zerosumfree semiring is also called an
\emph{antiring}~\cite{Tan2007}, which is thus another term that can be
used to describe $(\mathcal{P},+,\times)$.

\begin{remark} \label{remark2}

  $(\mathcal{P},+,\times)$ is \emph{entire}, because there are no
  non-zero elements \(p, q \in \mathcal{P}\) such that \(p \times q =
  \top\).

\end{remark}

This likewise shows that the monoid \((\mathcal{P},\times)\) is
completely removed from being a group, as there is no non-trivial
multiplicative inverse.

\begin{remark} \label{remark3}

  $(\mathcal{P},+,\times)$ is \emph{simple}, because \(\bot\) is
  \emph{infinite}, i.e., \(p + \bot = \bot, \forall p \in
  \mathcal{P}\).

\end{remark}


We may state another important definition~\cite{Golan1999} about
semirings, and observe a property of \((\mathcal{P},+,\times)\).

\begin{definition} \label{def_center}

  The \emph{center} \(C(\mathcal{P})\) of $\mathcal{P}$ is the set
  \(\{p \in \mathcal{P} \ | \ p \times q = q \times p, \ \mathrm{for \
    all} \ q \in \mathcal{P}\}\).

\end{definition}

\begin{remark} \label{remark4}

  The semiring $(\mathcal{P},+,\times)$ is \emph{commutative} because
  \(C(\mathcal{P}) = \mathcal{P}\).

\end{remark}

\section{Partial Ordering on a Boolean Propositional
  Algebra} \label{partialordering}

Consider a partial ordering relation \(\preccurlyeq\) on
\(\mathcal{P}\).  Informally, \(p \preccurlyeq q\) means that $p$ has
a lower measure of some metric than $q$ (e.g., $p$ is less likely to
be true than $q$, or is a weaker proposition than $q$).

Formally, $\preccurlyeq$ is a partial ordering on the semiring
$(\mathcal{P},+,\times)$ where the following conditions are
satisfied~\cite{HebWein1998}.

\begin{definition} \label{defpartialordering}

  If $(\mathcal{P},+,\times)$ is a semiring and \((\mathcal{P},
  \preccurlyeq)\) is a poset, then
  \((\mathcal{P},+,\times,\preccurlyeq)\) is a \emph{partially ordered
    semiring} if the following conditions are satisfied for all \(p,
  q,\) and $r$ in $\mathcal{P}$.

\begin{itemize}

\item[(i)] The \emph{monotony law of addition}:

\[p \preccurlyeq q \longrightarrow p + r \preccurlyeq q + r\]

\item[(ii)] The \emph{monotony law of multiplication}:

\[p \preccurlyeq q \longrightarrow p \times r \preccurlyeq q \times r.\]

\end{itemize}

\end{definition}

It is assumed that \(\top \preccurlyeq p, \forall p \in
\mathcal{P}\), and that \(p \preccurlyeq \bot\).

A semiring with a partial order defined on it is denoted as
\((\mathcal{P}, +, \times, \preccurlyeq)\).

Given Definition~\ref{defpartialordering}, it is instructive to
consider the behavior of the partial order under composition.  We
begin with a couple of simple results.

\begin{lemma} \label{lemma1}

  Given a partially-ordered semiring \((\mathcal{P}, +, \times,
  \preccurlyeq)\), \(\forall p, q \in \mathcal{P}\):

\begin{itemize}

\item[(i)] \(p \preccurlyeq p + q\), and
\item[(ii)] \(p \times q \preccurlyeq q\).

\end{itemize}

\end{lemma}

\begin{proof}

  For (i), consider that \(\top \preccurlyeq q\).  Using the monotony
  law of addition, we get \(\top + p \preccurlyeq q + p\).
  Considering that $\top$ is the additive identity element and that
  addition is commutative, we get \(p \preccurlyeq p + q\).

  For (ii), consider that \(q \preccurlyeq \bot\).  Using the monotony
  law of multiplication, we get \(q \times p \preccurlyeq \bot \times
  p\).  Considering that $\bot$ is the multiplicative identity element
  and that multiplication is commutative, we get \(p \times q
  \preccurlyeq p\). \qedhere

\end{proof}

Given these, we can state the following result on \((\mathcal{P}, +,
\times, \preccurlyeq)\).

\begin{theorem} \label{genconsistency}

Given \(p, q, r \in \mathcal{P}\),

\begin{itemize}

\item[(i)] if \(p + q \preccurlyeq r\), then \(p \preccurlyeq r\) and \(q
  \preccurlyeq r\); and

\item[(ii)] if \(p \preccurlyeq q \times r\), then \(p \preccurlyeq q\)
  and \(p \preccurlyeq r\).

\end{itemize}

\end{theorem}

\begin{proof}

  For part (i): The proof is by contradiction.  Assume the contrary.
  Then \(p + q \preccurlyeq r\), and at least one of \(p \preccurlyeq
  r\) or \(q \preccurlyeq r\) is false.

  Without loss of generality, assume that \(r \preccurlyeq p\).  Using
  the monotony law of addition and the commutativity of the $+$
  operator, \(q + r \preccurlyeq p + q\).

  Now, by Lemma~\ref{lemma1} (i), \(r \preccurlyeq q + r\).  Given
  the transitivity of $\preccurlyeq$, we get \(r \preccurlyeq p + q\),
  which is a contradiction.

  For part (ii): The proof is again by contradiction.  Assume the
  contrary.  Then \(p \preccurlyeq q \times r\) and at least one of
  \(p \preccurlyeq q\) and \(p \preccurlyeq r\) is false.

  Without loss of generality, assume that \(q \preccurlyeq p\).  Using
  the monotony law of multiplication and the commutativity of the
  $\times$ operator, we get \(q \times r \preccurlyeq p \times r\).
  
  Now, by Lemma~\ref{lemma1} (ii), \(p \times r \preccurlyeq p\).
  Given the transitivity of $\preccurlyeq$, we get \(q \times r
  \preccurlyeq p\), which is a contradiction. \qedhere

\end{proof}

The following result is similar.

\begin{theorem} \label{theorempartialorder}

  Given \(p, q, r, s \in \mathcal{P}\), if \(p \preccurlyeq q\) and
  \(r \preccurlyeq s\), then,

\begin{itemize}

\item[(i)] \(p + r \preccurlyeq q + s\), and
\item[(ii)] \(p \times r \preccurlyeq q \times s\).

\end{itemize}

\end{theorem}

\begin{proof}

These results can be proven directly.  Only (i) is proved, the proof
of (ii) being very similar.

We know the following:
\begin{eqnarray} \label{pleqq}
p \preccurlyeq q
\end{eqnarray}

and:
\begin{eqnarray} \label{rleqs}
r \preccurlyeq s
\end{eqnarray}

From~(\ref{pleqq}) and the monotony law of addition (considering the
direct sum of $s$ and both sides), we have:
\begin{eqnarray} \label{pq3}
p + s \preccurlyeq q + s.
\end{eqnarray}

Similarly, from~(\ref{rleqs}) and the monotony law (considering the
direct sum of $p$ and both sides), we have:
\begin{eqnarray} \label{rs4}
p + r \preccurlyeq p + s.
\end{eqnarray}

By considering transitivity in respect of (\ref{rs4}) and (\ref{pq3}),
we get \(p + r \preccurlyeq q + s\). \qedhere


\end{proof}

\begin{remark} \label{positivecone}

  The \emph{positive cone} $\hat{P}$ of \((\mathcal{P}, +,
  \preccurlyeq)\), which is the set of elements \(p \in \mathcal{P}\)
  for which \(p \preccurlyeq p + q, \forall q \in \mathcal{P}\), is
  the set $\mathcal{P}$ itself.  The \emph{negative cone} is empty.

\end{remark}

This is a direct consequence of Lemma~\ref{lemma1} (i), and it also
follows that the set of elements \(\{p \, | \, p + q \preccurlyeq p\}
= \emptyset\), showing that the negative cone is empty.

The analogous property of $\mathcal{P}$ in consideration of the
$\times$ operator can also be noted.

\begin{definition} \label{RTA}

\begin{itemize}

\item[(i)] A semiring $(\mathcal{P},+,\times)$ is called
  \emph{additively cancellative} if $(\mathcal{P},+)$ is cancellative,
  i.e., if $a + x = a + y$ implies $x = y$ for all $a,x,y \in
  \mathcal{P}$

\item[(ii)] Let $(\mathcal{P},+,\times)$ be a semiring with a zero
  $\top$. Then $\top$ is called \emph{multiplicatively absorbing} if
  $\top$ is absorbing in $(\mathcal{P},\times)$, i.e., if $\top a = a
  \top = \top$ holds for all $a \in \mathcal{P}$.

\end{itemize}

\end{definition}

We have the following result.

\begin{theorem} \label{theorempartialorder3}

  Let \(\mathcal{B} = (\mathcal{P}, +, \times,
  \preccurlyeq_{\mathcal{B}})\) be a partially ordered BPA, and
  $\widehat{\mathcal{B}} = (\widehat{\mathcal{P}}, +, \times)$ where
  $\widehat{\mathcal{P}} \subseteq \mathcal{P}$, is a BPSA of
  \((\mathcal{P}, +, \times)\).  Then,
	
\begin{itemize}

\item[(i)] if $(\mathcal{P}, +, \times)$ contains an additively
  cancellable element, but not $\top$, and
  
\item[(ii)] $\widehat{\mathcal{B}} \cup \left\{\top\right\} =
  \mathcal{B}$ if $(\mathcal{P}, +, \times)$ has a $\top$,
	
\end{itemize}
  
  then $\widehat{\mathcal{B}} = \mathcal{B}$.
  
\end{theorem}

\begin{proof}

  First we show that $\widehat{\mathcal{B}}$ is partially ordered.
  	
\begin{itemize}
	
\item[(a)] $p \preccurlyeq p$ (reflexivity) \\
  $p \preccurlyeq p$ always holds good, for some $p \in
  \widehat{\mathcal{B}}$
	
\item[(b)] if $p \preccurlyeq q$ and $q \preccurlyeq p$ then $p = q$
  (antisymmetry). \\ Since $p \preccurlyeq q$, therefore by
  Lemma~\ref{lemma1} (ii), $p + q = p$. Similarly $q + p = q$ Hence $p
  = q$ for some $p,q \in \widehat{\mathcal{B}}$.
	
\item[(c)] if $p \preccurlyeq q$ and $q \preccurlyeq r$ then $p
  \preccurlyeq r$ (transitivity) for some $p,q,r \in
  \widehat{\mathcal{B}}$. \\ By associative property of proposition
  and by Theorem~\ref{theorempartialorder}, it is transitive.
	
\end{itemize}
	
Hence, \((\widehat{\mathcal{P}}, +, \times,
\preccurlyeq_{\widehat{\mathcal{B}}})\) is partially ordered BPSA. So,
by Lemma~(\ref{lemma1}) we can say that $p \preccurlyeq q$ implies $p
+ x = y$ for some $x \in \widehat{\mathcal{B}}$, which satisfies the
monotony law of addition.  Also, if $p < q$ implies $p \times r
\preccurlyeq q \times r$ and $r \times p \preccurlyeq r \times q$ for
all $p,q \in \mathcal{B}$ and $r \in \widehat{\mathcal{B}}$ which in
turn satisfies monotony law of multiplication.  Therefore, we can say
that $\widehat{\mathcal{B}} \subseteq \mathcal{B}$ which means
$\preccurlyeq_{\mathcal{B}}$ and
$\preccurlyeq_{\widehat{\mathcal{B}}}$ similar.
	
If \((\mathcal{P}, +, \times)\) has a cancellable element, say $\top$,
then $\mathcal{B} \, \emptyset \, \widehat{\mathcal{B}}$ is either
empty or contains a single element, which has to be the $\top$ of
\((\mathcal{P}, +, \times)\). Also $\mathcal{B} \, \emptyset \,
\widehat{\mathcal{B}} = \phi$ and therefore $\widehat{\mathcal{B}} =
\mathcal{B}$.
		
We have $\top$ as a cancellable element of \((\mathcal{P}, +,
\times)\) and obtain $\mathcal{B} \, \emptyset \, \widehat{\mathcal{B}}
\subseteq \left\{\top\right\}$ , i.e., $\widehat{\mathcal{B}} \cup
\left\{\top\right\}$ \qedhere

\end{proof}

\section{Homomorphism and Isomorphism on Boolean Propositional
  Algebra} \label{congruences}

Let $\mathcal{Q} = (\mathcal{P},+,\times)$ and $\mathcal{R} =
(\mathcal{P},\oplus,\otimes)$ be two BPAs, then a mapping $\psi :
\mathcal{Q} \rightarrow \mathcal{R} \ \mathrm{of} \ \mathcal{Q}$ into
$\mathcal{R}$ is called a homomorphism of $(\mathcal{P},+,\times)$
into $(\mathcal{P},\oplus,\otimes)$ if,

\begin{itemize}

\item[(i) ] $\psi(a + b) = \psi(a) + \psi(b)$ and
\item[(ii) ] $\psi(a \times b) = \psi(a) \times \psi(b)$
are satisfied for all $a,b \in \mathcal{P}$. 
\end{itemize}

In other words we can say that $\psi$ is:

\begin{itemize}

\item[(i)] Order preserving: For each $x,y \in \mathcal{P}$, if $x
  \preccurlyeq y$, then $\psi(x) \preccurlyeq \psi(y)$

\item[(ii)] Operator preserving: For some operator $o$
  and each $x_{1},\ldots,x_{n} \in \mathcal{P}$,
  $\psi(o(x_{1},\ldots,x_{n})) = o(\psi(x_{1}),\ldots,f(x_{n}))$
  
\item[(iii)] Each mapping $\psi$ from $(\mathcal{P},+,\times)$ into
  $(T,\oplus,\otimes)$ determines an equivalence relation $\tau$ on
  $\mathcal{P}$ by $\tau = \psi^{-1} \odot \psi$, which may also be
  expressed by \(a \tau \acute{a} \equiv \psi(a) = \psi(\acute{a})\) for
  all \(a,\acute{a} \in \mathcal{P}\).

\end{itemize}

\begin{definition} \label{isomorphism} 

\begin{itemize}

\item[(i) ] An \emph{isomorphism} $\psi: \mathcal{Q} \rightarrow
  \mathcal{R}$ is a homomorphism such that the inverse map $\psi^{-1}:
  \mathcal{R} \rightarrow \mathcal{Q}$ - given by setting
  $\psi^{-1}(y) = x$ where $\psi(x) = y$ is a homomorphism. Two BPAs
  are isomorphic if and only if there is an isomorphism from one to
  the other.

\item[(ii) ]  Let \((\mathcal{P},+,\times)\) be BPA and \(\psi:
  (\mathcal{P},+,\times) \rightarrow (\mathcal{P},\oplus,\otimes)\) is
  a homomorphism.  Then \((\psi(\mathcal{P}),+,\times)\) is again a BPA.
\end{itemize}

\end{definition}

Based on this, we have the following.

\begin{theorem} \label{theoremhomorphism2}

  Let
  $(\mathcal{P},+,\times),(\mathcal{P}_{1},+,\times),(\mathcal{P}_{2},+,\times)$
  are BPAs and there exists two homomorphic functions such that,
  \(\psi_{1}: \mathcal{P} \rightarrow \mathcal{P}_{1}\) and
  \(\psi_{2}: \mathcal{P} \rightarrow \mathcal{P}_{2}\), then the
  homomorphic function \(\psi: \mathcal{P}_{1} \rightarrow
  \mathcal{P}_{2}\) establish the relation of BPA and BPSA between
  $\mathcal{P}_{1},\mathcal{P}_{2}$.
  
  \end{theorem}

\begin{proof}

  The proof is by contradiction.  Assume to the contrary that
  homomorphism $\psi$ exists.  Since $\psi_{1}$ is surjective, for
  $a_{1} \in \mathcal{P}_{1}$, there is some $a \in \mathcal{P}$ which
  satisfies $\psi_{1}(a) = a_{1}$. Also we know that $\psi \odot
  \psi_{1} = \psi_{2}$ by Definition \ref{isomorphism}(ii), we can
  show that $\psi(a_{1}) = \psi(\psi_{1}(a)) = \psi_{2}(a)$.  Hence,
\begin{eqnarray}\label{heq1}
\psi(a_{1}) = \psi_{2}(a) \ \forall a \ \in \mathcal{P} \ \mathrm{such \ that} \ \psi_{1}(a) = a_{1}.
\end{eqnarray}
This shows that $\psi_{1}(a) = \psi_{1}(\acute{a}) \implies
\psi_{2}(a) = \psi_{2}(\acute{a})$ where $a,\acute{a} \in
\mathcal{P}$. Therefore, \(\tau_{1} \subseteq \tau_{2}\)

We have to show $\psi$ is surjecive iff $\psi_{2}$ is surjective.

To prove by contradiction, assume that on the contrary we have
\(\tau_{1} \subseteq \tau_{2}\).

For $\psi$ to be surjective, necessary conditions are:

\begin{itemize}

\item[(i)] $\psi(a_{1} + b_{1}) = \psi(a_{1}) + \psi(b_{1})$ and
\item[(ii)] $\psi(\top) = \top$
\item[(iii)] $\psi(\bot) = \bot$
are satisfied for all $a_{1},b_{1} \in \mathcal{P}$

\end{itemize}

By (\ref{heq1}), we have $\psi_{1}(a) = a_{1} = \psi_{1}(\acute{a})
\implies a \tau_{1} \acute{a} \implies a \tau_{2} \acute{a}$ This in
return shows that $\psi_{2}(a) = \psi_{2}(\acute{a})$

Hence, $\psi$ defines the mapping of $\mathcal{S}_{1}$ into
$\mathcal{S}_{2}$ and $\psi_{2}(a) = \psi(\psi_{1}(a)) \implies \psi
\odot \psi_{1} = \psi_{2}$

Also, $\psi(\top + a) = \psi(\top) + \psi(a)$ for some $a \in
\mathcal{P}$ Since, $\psi(\top) \preccurlyeq \psi(a)$. \qedhere

\end{proof}

\begin{theorem} \label{theoremhomorphism}

  Let $(\mathcal{P},+,\times)$ and $(\mathcal{P}_{1},+,\times)$ are
  two BPAs and there exist a homomorphic function \(\psi: \mathcal{P}
  \rightarrow \mathcal{P}_{1}\). The function \(\psi\) is an onto
  order preserving iff it is an isomorphism.
  
\end{theorem}
\begin{proof}

  If $\psi$ is an onto order preserving, we first need to show that
  $\psi^{-1}$ is well defined.  Since $\psi$ is onto, for every $y$
  there is at least one $x \in \mathcal{P}$ where $\psi(x)=y$.  Since
  $\psi$ is an order preserving, it is one to one, so there is not
  more than one $x$ where $\psi(x)=y$. Hence by definition of
  $\psi^{-1}$ we can show that $\psi(\psi^{-1}(y))=y$ and
  $\psi^{-1}(\psi(x))=x$ for each $y \in \mathcal{P}_{1}$ and $x \in
  \mathcal{P}$.

  Now $\psi^{-1}$ is order preserving since if $\psi^{-1}(y) \notin
  \psi^{-1}(\acute{y})$, we have $\psi(\psi^{-1}(y)) \notin
  \psi(\psi^{-1}(\acute{y}))$ implies $y \notin \acute{y}$.
  Contraposing this, we have $y \preccurlyeq \acute{y}$ only if
  $\psi^{-1}(y) \preccurlyeq \psi^{-1}(\acute{y}$.

  Similarly, $\psi^{-1}$ preserve operators. For any $n$ place
  operator $\tau$, $\psi^{-1}(\tau(y_{1},\ldots,y_{n})) \notin
  \tau(\psi^{-1}(y_{1}),\ldots,\psi^{-1}(y_{n})$, then since $\psi$ is
  one to one, we have $\psi(\psi^{-1}(\tau(y_{1},\ldots,y_{n})))
  \notin \tau(\psi(\psi^{-1}(y_{1})),\ldots,\psi(\psi^{-1}(y_{n})))$.

  But, we have $\psi(\psi^{-1}(y)) = y$ for all $y$, we get
  $\tau(y_{1},\ldots,y_{n}) \notin \tau(y_{1},\ldots,y_{n})$, which is
  contradiction Therefore, $\psi^{-1}$ must preserve operators, so
  $\psi$ is an isomorphism.

  Conversely, if $\psi$ is an isomorphism, then we have to show that
  it is an onto order preserving.  Primarily, If $y \in
  \mathcal{P}_{1}$, we have $\psi(\psi^{-1}(y))=y$, and hence $\psi$
  is onto.  Secondly, if $\psi(x) \preccurlyeq \psi(y)$, then
  $\psi^{-1}(\psi(x)) \preccurlyeq \psi^{-1}(\psi(y))$, will give $x
  \preccurlyeq y$, and hence, $\psi$ is an order preserving. \qedhere

\end{proof}

\section{Differences in Boolean Propositional Algebra}

Given any BPA $\mathcal{B} = (\mathcal{P},+,\times,\preccurlyeq)$, its
\emph{ideal BPA}~\cite{Henrik} \cite{Iseki} is $Ideal(\mathcal{B})$,
ordered under $\subseteq$ and with operators $(+,\times)$.

Let $\mathcal{B} = (\mathcal{P},+,\times,\preccurlyeq)$ be BPA and
(subtrahends $\ominus$) ~\cite{Dirk} be BPSA of $\mathcal{B}$ whose
elements are cancellable in $(\mathcal{P},+,\times,\preccurlyeq)$.  We
further assume the existence of a zero $\top$ and contains an opposite
element $\neg \alpha$ for each $\alpha \in \ominus$. Then without
restriction of generality, $\ominus$ can be chosen as
$Ideal(\mathcal{B})$ of $(\mathcal{P},+,\times)$. For convinience, we
can write BPA of differences w.r.t $\mathcal{P}$ as
$D(\mathcal{P},\ominus)$
 
\begin{definition} \label{facts2}

  Let $(\mathcal{P},+,\times,\preccurlyeq)$ be a BPA and
  $(\mathcal{P},\oplus,\otimes) = D(\mathcal{P},\ominus)$ a BPA of
  differences of $(\mathcal{P},+,\times)$, then $p \preccurlyeq q
  \mathrm{implies \ there \ exists \ some} \ \Delta \in \ominus \
  \mathrm{and} \ p,q \in \mathcal{P} \ \mathrm{such \ that}$
\begin{eqnarray} \label{propdiffeqn1}
p + \Delta \preccurlyeq q + \Delta
\end{eqnarray}
defines the smallest extension of $\preccurlyeq$ within
$\mathcal{P}$. For convinience we can denote this smallest extension
as $\preccurlyeq^{\prime}$. Further
$(\mathcal{P},+,\times,\preccurlyeq^{\prime})$ is a partial order BPA
with the property
\begin{eqnarray} \label{propdiffeqn2}
p \preccurlyeq^{\prime} q \Leftrightarrow p + \xi \preccurlyeq^{\prime} q + \xi \ \mathrm{for \ all} \ p,q 
\in \mathcal{P} \ and \ \xi \in \ominus
\end{eqnarray}
Moreover, $\preccurlyeq^{\prime}$ and $\preccurlyeq$ are similar iff
$(\mathcal{P},+,\times,\preccurlyeq)$ itself satisfies
\begin{eqnarray} \label{propdiffeqn3}
p \preccurlyeq q \Leftrightarrow p + \xi \preccurlyeq q + \xi \ \mathrm{for \ all} \ p,q \in \mathcal{P} \ \mathrm{and} \
 \xi \in \ominus
\end{eqnarray}

\end{definition}

We have the following result.

\begin{theorem} \label{propdiff}

  Let $(\mathcal{P},\oplus,\otimes) = D(\mathcal{P},\ominus)$ be BPA
  of differences of BPA $(\mathcal{P},+,\times)$ with respect to the
  ideal of subtrahends $\ominus$, then if $(\mathcal{P},+,\times)$ is
  multiplicatively left cancellative, the necessary and sufficient
  condition for $D(\mathcal{P},\ominus)$ to be multiplicatively
  left-cancellative is
\begin{eqnarray} \label{propdiffeqn}
\Delta \neq c \ and \ a \neq b \ \Rightarrow \ c \times a + \Delta b \neq c \times b + \Delta a
\end{eqnarray}
for all $a,b,c \in \mathcal{P}$ and $\Delta \in \ominus$
\end{theorem}

\begin{proof}

  Let $(\mathcal{P},\oplus,\otimes)$ be multiplicatively left
  cancellative. Then $\Delta \neq c$ and $a \neq b$ imply $(c -
  \Delta)a \neq (c - \Delta)b$ and thus $c \times a + 
  \Delta b \neq c \times b + \Delta a$, which 
  proves (\ref{propdiffeqn}). For the converse we assume $(c - \Delta)
  (a - \alpha) = (c - \Delta)(b - \alpha)$ for arbitrary elements 
  $a - \alpha, b - \alpha$ and $c - \Delta \neq \top$ of $\mathcal{P}$. 
  This yields $c \times a + \Delta \alpha + c \alpha +
  \Delta b = c \alpha + \Delta a + c \times b + 
  \Delta \alpha$.  Since $\Delta \alpha \in \ominus$ and $c 
  \alpha \in \ominus$ are cancellable in
  $(\mathcal{P},+)$, we obtain $c \times a + \Delta b = c \times b + \Delta a$. This
  and $c \neq \Delta$ imply $a = b$ by (\ref{propdiffeqn}) and hence
  $a - \alpha = b - \alpha$. Therefore (\ref{propdiffeqn}) implies
  that $(\mathcal{P},\oplus,\otimes) = D(\mathcal{P},\ominus)$ is
  multiplicatively left cancellative, which yields the same for
  $(\mathcal{P},+,\times)$ and completes the proof. \qedhere

\end{proof}

\bibliography{paper}

\end{document}